\documentclass{cmslatex}
\usepackage[paperwidth=7in, paperheight=10in, margin=.875in]{geometry}
 \usepackage[backref,colorlinks,linkcolor=red,anchorcolor=green,citecolor=blue]{hyperref}
\usepackage{amsfonts,amssymb}
\usepackage{amsmath}
\usepackage{bm}
\usepackage{mathrsfs,esint}
\usepackage{graphicx}
\usepackage{cite}
\usepackage{enumerate}
\usepackage{tikz}
\usetikzlibrary{decorations.pathreplacing, decorations.pathmorphing, decorations.shapes}
\usepackage{algorithm2e}
\renewcommand{\theequation}{\arabic{section}.\arabic{equation}}

\newcommand{\beq}{\begin{equation}}
\newcommand{\eeq}{\end{equation}}

\newcommand{\ga}{{\gamma}}

\newcommand{\del}{\delta}

\newcommand{\Om}{\Omega}
\newcommand{\om}{\omega}

\def\cH{\mathcal{H}}

\def\cL{\mathcal{L}}

\def\eb{\mathbf{e}}
\def\ub{\mathbf{u}}

\def\fb{\mathbf{f}}

\def\bx{\bm x}
\def\by{\bm y}

\def\R{\mathbb{R}}

   \allowdisplaybreaks
\begin{document}
 \title{Fast algorithm for computing nonlocal operators with finite interaction distance\thanks{Received date, and accepted date.}}


         \author{Xiaochuan Tian\thanks{Department of Mathematics, The University of Texas at Austin, Austin, TX 78712 (xtian@math.utexas.edu).} \and Bj\"{o}rn Engquist\thanks{ Department of Mathematics and the Oden Institute, University of Texas at Austin, Austin, TX 78712 (engquist@ices.utexas.edu).}
}

         \pagestyle{myheadings} \markboth{FAST ALGORITHM FOR NONLOCAL OPERATORS}{XIAOCHUAN TIAN AND BJ\"{O}RN ENGQUIST }
         
          \maketitle

\begin{abstract}
              Developments of nonlocal operators for modeling processes that traditionally have been described by local differential operators have been increasingly active during the last few years. One example is peridynamics for brittle materials and another is 
nonstandard diffusion including the use of fractional derivatives. A major obstacle for application of these methods is the high computational cost from the numerical implementation of the nonlocal operators. It is natural to consider fast methods of fast multipole or hierarchical matrix type to overcome this challenge. Unfortunately the relevant kernels do not satisfy the standard necessary conditions. In this work a new class of fast algorithms is developed and analyzed, which is some cases reduces the computational complexity of applying nonlocal operators to essentially the same order of magnitude as the complexity of standard local numerical methods.
 \end{abstract}

          \section{Introduction}\label{intro}
          Nonlocal operators as alternatives to the local differential operators have been widely applied to the modeling of a 
          number of physical and biological processes over the last 
          few years.  One example is given by the peridynamics models for brittle materials 
          \cite{Silling2000,Bobaru2008,Foster2010,Bobaru2010a,Bobaru2011a,Silling2007}, and another 
          is nonstandard diffusion including the use of fractional derivatives \cite{Andreu2010,Bates1999,Du2012a,Fife2003,MeKl00,Valdinoci2017}. 
          Moreover, even in the case of solving differential equations, nonlocal integral relaxations to local differentiations are sometimes introduced 
as a numerical technique, see \cite{Gingold1977,Shi2016,Nochetto2018b,Nochetto2018a}. 

While the nonlocal models can provide more accurate descriptions of a physical system, the nonlocality
also increases the computational cost compared to classical local models
based on partial differential equations (PDEs).
Thus it is imperative to develop fast algorithms for the computation of nonlocal models. 
In this work, we focus on nonlocal models with finite nonlocal interaction distances (see \cite{Tian2016a}).
Each point in the domain interacts with points within certain ``horizon'' of it. 
That is to say the interaction kernel is compactly supported in certain finite region and may also be singular at zero distance. 
When the interaction region is spherical, we call the radius of the region the horizon radius or just horizon in short.
Such kernels make most fast algorithms that utilizing the low rank property of the far field interactions fail to work, because 
the kernel is typically not regular around the boundary of its support. 
In the work of \cite{WangHong2016,WangHong2017}, FFT based fast algorithms are developed for solving peridynamic models,
in which Toeplitz structures of matrices are greatly exploited. 
However, such methods fail to work for the more generals cases of inhomogeneuous media and variable horizons
resulting in non-translation invariant kernels, which is the concern of this work. 
 
The peridynamics (PD) theory of solid mechanics takes into account of large scale deformations and 
reformulates the internal force of a deformed body by introducing the nonlocal way of interactions of materials points. 
Let $\ub$ denote the displacement field, 
the force balance equation in the bond-based peridynamic theory is given by 
\beq\label{eq:PDmain}
- \int_{\cH_{\bx}} \fb (\ub(\by)-\ub(\bx), \by, \bx) d\by= \mathbf{b}  (\bx)\,, 
\eeq
where the integration is limited to the finite region $\cH_{\bx}$.
$ \mathbf{b}= \mathbf{b} (\bx)$ denotes the body force and $\fb$ is the force density function, which contains all the constitutive properties.
For a prototype micro-elastic material,  the force vector $\fb$ points to the direction of the deformed bond 
with its magnitude a function of the bond relative elongation(stretch). In general, $\fb$ depends nonlinearly on 
the bond relative elongation, which gives the model the flexibility to describe changes of material properties. In this work, we will only consider the linear model obtained from the assumption of small deformation, namely, 
\beq
 \fb (\ub(\by)-\ub(\bx), \by, \bx) =  \om(\bx, \by) [\eb\cdot(\ub(\by)-\ub(\bx))]\eb\,,
\eeq
where $\eb$ is the unit vector in the bond direction $\by-\bx$, and $ \om(\bx, \by)$ is the micromodulus function satisfying 
\[
\int_{\cH_{\bx}} |\by-\bx|^2\om(\bx, \by) d\by< \infty \,.
\]
For each $\bx$,  $\om(\bx, \cdot)$ is a function with support on $\cH_{\bx}$. 
In practice, this is achieved by choosing a fractional type kernel multiplied
with the characteristic function or the conical function on $\cH_{\bx}$ (see \cite{Bobaru2010a}), where
both choices create non-smooth transition to zero at the boundary $\partial \cH_{\bx}$.

As we can see, the discretization of the nonlocal operator in \eqref{eq:PDmain} results in a matrix 
$A$ with high density, for which fast algorithms must be considered in order to lower the cost of matrix multiplication and inversion. 
There are several classical algorithms for the treatment of non-sparse matrices, including the 
fast multipole methods (FMMs) \cite{GrRo87,GrRo97,YBZ04}
and the hierarchical matrix techniques \cite{Hackbusch99,Hackbusch02}. 
More recent developments can be found in \cite{HoYi16}. 
We will first consider the existing fast solvers to see whether they are effective in the computation of nonlocal models like \eqref{eq:PDmain}. 
In Section \ref{sec:2} we will show that the performance of the fast solvers depends largely on the
regularity of the kernel across the boundary $\partial \cH_{\bx}$.
This sends the clear message of the preference of the far field behavior of the kernel function in terms of computation efficiency. 
Popular choices of kernels in the peridynamics simulation will cause serious problems in the interest of fast computation. 
We then propose in Section \ref{sec:3} to split the two types of singularities -- one at origin and the other at the boundary $\partial \cH_{\bx}$-- of the kernel function, 
and treat them separately. 
Once the kernel is decomposed into two, one that is smooth away from origin and the other that is smooth away from the truncation, 
we then use different fast algorithms that exploit the different types of smoothness of the two kernels. 
The classical fast algorithms such as FMMs and hierarchical matrix techniques only resolve the first type of singularity, namely that the singularity set is one point (the origin of the kernel in most cases). 
Notice that the second type of singularity is on the boundary set $\partial \cH_{\bx}$, and it is a set of codimension $1$. In 2d the sets of singularity are boundary curves, and in 3d they are boundary surfaces.
Our algorithm that deals with second type of singularity is basically a new FMM type method for kernels
that exhibit singularities on codimension 1 sets. There is a reduction of computational complexity in all dimension with optimal rate $O(N\log N)$ for $N$
unknowns in 1d. 
Other potential applications of this work include computing with the 
retarded potentials raised in time domain boundary
integral equations \cite{Sayas2016}, where the potentials are discontinuous functions defined in space-time.

\section{The role of far field behaviors of kernels} \label{sec:2}
We will use the nonlocal diffusion operator to illustrate our methodology in this section. The nonlocal diffusion operator is given by
\beq \label{eq:diffusion}
\cL u(\bx) = \int_{\cH_{\bx}}\om(\bx,\by) (u(\by)-u(\bx)) d\by \,.
\eeq 
See more discussions of details in \cite{Du2012a}. Here we assume that the kernel $\om$ is given by 
\beq \label{kernel}
\om(\bx, \by)=\frac{C(\bx,\by)}{\del(\bx,\by)^{d+2}} \gamma(\frac{|\by-\bx|}{\del(\bx,\by)})\,,
\eeq
where the kernel $\gamma$ is a nonnegative function supported on $[-1,1]$ with $\int_{-1}^1 s^2 \gamma(|s|) ds =1$.
$C(\bx,\by)$ is the diffusion coefficient and $\del(\bx,\by)$ is the horizon function that satisfy
\beq
C_0\leq C(\bx,\by)\leq C_1, \text{ and } 0<\del(\bx,\by)\leq \del_1 \,.
\eeq
Here we keep the dependence of $C$ and $\del$ on $\bx$, $\by$
so that it could model the inhomogeneous medium, see related works \cite{Silling2015a,Tian2016b,DuZZ}.  

The common practice in the peridynamics simulation (see \cite{Bobaru2010a}) takes $\ga$ to be
\beq \label{kernelgamma}
\gamma(|s|)=\frac{c_1}{|s|}\chi(|s|<1), \text{ or } \gamma(|s|)=\frac{c_2}{|s|}(1-|s|)\chi(|s|<1)\,. 
\eeq
We distinguish between the two cases: 
\begin{enumerate}[(a)]
\item If $C$ and $\del$ only depends on $x$, namely $C(\bx,\by)= C(\bx)$ and $\del(\bx,\by)= \del(\bx)$, then \eqref{eq:diffusion} is the nonlocal diffusion operator
of non-divergence type. As $\del_1\to0$, we have
\[
\cL u(\bx)  \to C(\bx) \Delta u (\bx)\,.
\]
\item If $C$ and $\del$  are symmetric, namely $C(\bx,\by)= C(\by, \bx)$ and $\del(\bx,\by)= \del(\by, \bx)$, then \eqref{eq:diffusion} is the nonlocal diffusion operator
of divergence type. As $\del_1\to0$, we have
\[
\cL u(\bx)  \to  \nabla (\sigma(\bx) \nabla u (\bx))\,, \text{ where }\sigma(\bx)=C(\bx,\bx)\,.
\]
\end{enumerate}

The FFT-based methods mentioned earlier are restrictive in the application to heterogeneous media. Here we seek 
algorithms that can be applied to the more general cases given by \eqref{kernel}.
 We use the recently proposed hierarchical interpolative factorization (HIF) \cite{HoYi16} method to factorize the dense matrix
$A$ obtained from discretizing the nonlocal operator $\cL$. The factorization of the matrix can
be used to rapidly apply both $A$ and $A^{-1}$. Therefore it can serve as a direct solver or be
stored and reused for iterative methods or time-dependent problems. 
The HIF takes advantage of the low-rank behavior of the off-diagonal entries of the dense matrix, 
so the far field behavior of the kernel $\om$ is the central factor on the performance of the the algorithm. 
We will illustrate in this section that the performance of the HIF algorithm improves greatly as the smoothness of
 the kernel away from zero enhances. All computations are performed in MATLAB R2014b on a single core of a 1.1GHz Intel Core M CPU on a 64-bit Mac laptop.

\subsection{Nonlocal diffusion in homogeneous media} \label{sec2:homo}
In this section, we solve the the nonlocal problem on the interval $\Om=(0,1)$ with $C(x,y)\equiv 1$ and $\del(x,y)\equiv\del$
with $\del$ being the horizon radius of the homogeneous medium. 
The matrix $A$ is obtained form the discretization of the nonlocal operator $\cL$, which is is based on the asymptotically compatible schemes in \cite{Tian2013a}, which ensures uniform discretization error independent of the horizon function $\del$. In Table \ref{table:1}, the storage and the time of matrix-vector multiplication are listed for the original matrix $A$ as well as its compression after using the HIF method.
The original matrix $A$ is stored using the sparse matrix representation with the 
MATLAB built-in matrix-vector multiplication applied to the test. The kernel $\gamma(s)$ has $O(1/|s|)$ singularity at $0$ and is assumed to be smooth expect at origin and $\pm1$
 ($\om$ has singularity at $0$ and $\pm\del$). The number of discretization nodes $N$ is fixed to be $2048$ in Table \ref{table:1}. 
The relative precision parameter for the HIF method is set to be $1E-8$.
 
 \begin{table}[htbp]
\[
 \begin{array}{|cc|cc|cc|}
\hline
 &  & \multicolumn{2}{|c|}{\text{HIF}} &   \multicolumn{2}{|c|}{\text{Original}}\\ \hline
\text{kernel} & \delta & \text{Memory (MB)}& T_{\text{multi}} \text{(sec)} &  \text{Memory (MB)}& T_{\text{multi}} \text{(sec)} \\ 
C^{-1} & 1 & 1.93 & 3.54E-3 & 67.13 & 8.49E-3\\
C^{0} & 1 & 1.93 & 3.52E-3 & 67.13 & 8.61E-3\\
C^{1} & 1& 1.93 &3.53E -3& 67.13 &8.21E-3\\ 
C^{2}& 1 & 1.93&3.43E-3& 67.13 &8.21E-3\\ 
C^{3}&1 &1.93&3.38E-3 &67.13&8.25E-3\\ \hline
C^{-1} & 1/4 & 33.58 &6.23E -3 & 29.40 &3.34E-3\\ 
C^{0} & 1/4 & 33.37 &6.34E -3 & 29.40 &3.36E-3\\ 
C^{1} &1/4 &4.44&4.34E-3&29.40 &3.39E-3\\ 
C^{2}&1/4&2.10&3.61E-3&29.40&3.43E-3\\ 
C^{3}&1/4 &1.99&3.43E-3 &29.40&3.41E-3\\ \hline
C^{-1} & h &2.11&3.30E-3&0.11&5.04E-5\\
C^{0} & h &2.11&3.21E-3&0.11&4.57E-5\\
C^{1} & h&2.11&3.27E-3&0.11&4.57E-5\\ 
C^{2}&h&2.11&3.21E-3&0.11&4.76E-5\\ 
C^{3}&h &2.11&3.31E-3 &0.11&4.66E-5\\ \hline
\end{array} 
\]

\caption{The storage and time of matrix-vector multiplication for the 1d nonlocal diffusion in homogenous medium. The domain $\Om=(0,1)$ is discretized using $N=2048$ points. The mesh size $h=1/2048$. The relative precision of HIF approximation is set to be $\epsilon=10^{-8}$. } \label{table:1}
\end{table}

The left column of Table \ref{table:1} lists the different types of regularity (at $\pm1$) of the kernel  $\gamma(s)$   and 
the different choices of the horizon parameter $\del$.
The discontinuous kernel (denoted as $C^{-1}$) is the one that is most often used in peridynamics simulations given by 
\[
\ga(|s|)=\frac{1}{|s|}\chi(|s|<1)\,.
\]
The more regular $C^k$ ($k=0,1,2,3$)  kernels are obtained by subtracting $\gamma$ with polynomials $p^{2k}$ ($k=0,1,2,3$) of degree $2k$. The details of constructing these polynomials will be given in Section 3. When $\del\geq1$, the matrix $A$ is fully dense and it does not see the tail of the kernel. 
As a result, the HIF works well for all the types of kernels. When $\del$ gets smaller,  the regularity of the kernel $\gamma(s)$ at $\pm1$
plays a decisive role on the performance of the HIF. 
We see that when  $\gamma(s)$  is discontinuous  at $\pm1$, no compression of the matrix $A$ is observed using the HIF.
The storage taken and computation time decrease with increasing regularity of  $\gamma(s)$  from discontinuous ($C^{-1}$) 
all the way to three times continuously differentiable ($C^3$).
In the case of $\del=1/4$, the critical improvement of the computational efficiency happens between $C^0$ kernel functions and $C^1$ kernel functions
as shown in Table  \ref{table:1}. 
We remark that the location of critical regularity changes with the precision parameter $\epsilon$.
When the HIF algorithm runs with higher precision, the higher critical regularity is observed;
and when the HIF algorithm runs with lower precision,  the lower critical regularity is observed.
Similar patterns are also observed in the experiments that will be discussed in Section \ref{sec2:hetero} and \ref{sec2:2d}. 
When $\del$ keeps decreasing, the effect of different regularity of the kernel  $\gamma(s)$  will diminish. 
The extreme case is that $\del$ equals the mesh size $h$ when the matrix $A$ is in fact tridiagonal. 
The case of $\delta=h$ is listed in Table \ref{table:1} to compare the overhead of using the HIF algorithm
with the sparse matrix.  
The overhead of using the HIF algorithm in this case  is due to the fact that  HIF does not build into itself 
 the ability to detect the horizon and ignore the zeros outside. It will be future work to explore this possibility depending on application.

To account for the complexity of the HIF, we provide the computation time for different discretization nodes $N$ in Table \ref{table:2}. 
The horizon parameter $\del$ in this table is fixed to be $1/4$.  
The scaling results for the choice of the $C^0$ kernel is compared with the choice of the $C^3$ kernel. 
The complexity of the matrix-vector multiplication for the original matrix $A$ without compression is approximately $O(N^2)$ in both cases.
On the other hand, the HIF code behaves quite differently for the $C^0$ kernel and the $C^3$ kernel. 
The observed complexity of HIF with the $C^0$ kernel  is still , while the complexity of HIF with the $C^3$ kernel is just around $O(N)$. 

The simulations show effectiveness of the HIF fast algorithm with the improving far field regularity of nonlocal interaction kernels. 
Finally, we remark that the cost of constructing the factorization is in general larger than the cost of matrix-vector multiplication. 
This is acceptable since one only need to construct the factorization once and reuse it
for an iterative method or in  time-dependent problems.

 \begin{table}[htbp]
\[
 \begin{array}{|cc|cc|cc|}
\hline
  & &   \multicolumn{2}{|c|}{\text{HIF}} &   \multicolumn{2}{|c|}{\text{Original}}\\ \hline
\text{kernel} &N  & \text{Memory (MB)}&T_{\text{multi}} \text{(sec)}&  \text{Memory (MB)}& T_{\text{multi}} \text{(sec)}\\ 
C^0 &512  &2.10&4.87E-4&1.85&1.68E-4\\
C^0 &1024 &8.36&1.56E-3&7.36&8.76E-4 \\ 
C^0 &2048&33.35&5.31E-3&29.40&2.96E-3\\ 
C^0 & 4096 &122.99&2.11E-2&117.52&1.39E-2\\ \hline 
C^3 &512  &0.41&1.30E-3&1.85&1.63E-4\\
C^3 &1024 &1.04&1.55E-3&7.36&8.95E-4\\ 
C^3 &2048&1.99&3.07E-3&29.40&3.35E-3\\ 
C^3 & 4096 &3.93&5.97E-3  &117.52&1.48E-2\\ \hline
\end{array} 
\]
\caption{The storage and time of matrix-vector multiplication for the 1d nonlocal diffusion in homogenous medium. $\del=1/4$. The relative precision of HIF approximation is set to be $\epsilon=10^{-8}$.} \label{table:2}
\end{table}

\subsection{Nonlocal diffusion in heterogeneous media} \label{sec2:hetero}
 We test the effectiveness of the HIF algorithm for nonlocal diffusion in heterogeneous media in this section. The kernel $\om$ is defined as \eqref{kernel}, where $C(x,y)\equiv1$ and $\del(x,y)=\del(x)$. Each point $x$ has its own horizon radius $\del(x)$.
 In this case, $\cL$ is a nonlocal diffusion operator of non-divergence type. 
The matrix $A$ obtained form the discretization of the nonlocal operator $\cL$ is now non-symmetric, in which case the
FFT-based fast algorithm will not work. 
In Table \ref{table:3}, the storage and the time of matrix-vector multiplication are listed for the original matrix $A$ as well as its compression after using the HIF algorithm. 

\begin{figure}[htbp] 
  \centering
\includegraphics[width=6cm]{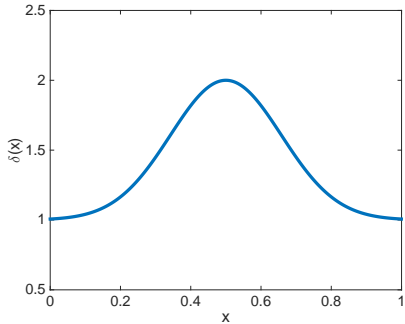}
\caption{Plot of the horizon function $\del(x)=\del_0(1+e^{-20(x-0.5)^2})$ with $\del_0=1$.} \label{fig:delta}
    \end{figure}

In  the numerical experiments we use $\del(x)=\del_0(1+e^{-20(x-0.5)^2})$ for $x\in \Om=(0,1)$ so that $\del(x)$ is a function between $\del_0$ and $2\del_0$, 
as shown in Fig \ref{fig:delta}. 

The kernel $\gamma(s)$ has $O(1/|s|)$ singularity at $0$ and is assumed to be smooth expect at origin and $\pm1$. 
The left column of Table \ref{table:3} lists the types of regularity (at $0$ and $\pm1$) of the kernel $\gamma$ and the horizon parameter $\del_0$. The HIF algorithm for nonlocal diffusion in heterogeneous media works similarly as the homogeneous cases.  In particular, we tailored Table \ref{table:3} to the case where $\del_0=1/4$ because for the extreme cases where $\del_0=1$ or $\del_0=h$, the variation of regularity of the kernel does not
result in any change in the computational cost, which is exactly the same behavior observed in Table \ref{table:1}.

 \begin{table}[htbp]
\[
 \begin{array}{|cc|cc|cc|}
\hline
 &  & \multicolumn{2}{|c|}{\text{HIF}} &   \multicolumn{2}{|c|}{\text{Original}}\\ \hline
\text{kernel} & \delta_0 & \text{Memory (MB)}& T_{\text{multi}} \text{(sec)} &  \text{Memory (MB)}& T_{\text{multi}} \text{(sec)} \\ 
C^{-1} & 1/4 & 55.90 &6.18E -3 & 40.24 &4.80E-3\\ 
C^{0} & 1/4 & 55.80 &6.24E -3 & 40.24 &4.72E-3\\ 
C^{1} &1/4 &5.10&3.92E-3&40.24 &4.77E-3\\ 
C^{2}&1/4&2.44&3.27E-3&40.24&4.79E-3\\ 
C^{3}&1/4 &2.31&3.13E-3 &40.24&4.67E-3\\ \hline
\end{array} 
\]
\caption{The storage and time of matrix-vector multiplication for the 1d nonlocal diffusion in heterogeneous medium. The domain $\Om=(0,1)$ is discretized using $N=2048$ points. The mesh size $h=1/2048$. The relative precision of HIF approximation is set to be $\epsilon=10^{-8}$. } \label{table:3}
\end{table}

To account for the complexity, we provide the computation time for different numbers of discretization nodes $N$ in Table \ref{table:4}. 
The horizon parameter $\del_0$ is fixed to be $1/4$.
Two types of kernels -- the $C^0$ kernel and the $C^3$ kernel -- are used 
in the simulation. 
The complexity of the matrix-vector multiplication for the original matrix $A$ without compression is $O(N^2)$ as expected. 
Similar to the homogeneous case, the observed complexity of HIF in the heterogeneous case is again $O(N^2)$ for the $C^0$ kernel and $O(N)$ for the $C^3$ kernel. The numerical examples in Sections \ref{sec2:homo} and \ref{sec2:hetero} 
indicate that the HIF has the same effectiveness for models in heterogeneous media  
as it has for models in homogeneous media if the kernel functions are sufficiently smooth away from origin.

 \begin{table}[htbp]
\[
 \begin{array}{|cc|cc|cc|}
\hline
 &&   \multicolumn{2}{|c|}{\text{HIF}} &   \multicolumn{2}{|c|}{\text{Original}}\\ \hline
\text{kernel}& N  & \text{Memory (MB)}&T_{\text{multi}} \text{(sec)}&  \text{Memory (MB)}& T_{\text{multi}} \text{(sec)}\\ 
C^0 &512  &3.50&4.48E-4&2.52&2.29E-4\\
C^0 &1024 &13.96&1.75E-3&10.07&1.12E-3\\ 
C^0 &2048&55.80&5.56E-3&40.24&4.13E-3\\ 
C^0 & 4096 &221.35&1.98E-2&160.87&1.61E-2\\ \hline 
C^3 & 512  &0.53&1.34E-3&2.52&2.42E-4\\
C^3 &1024 &1.22&1.42E-3&10.07&1.12E-3\\ 
C^3 &2048&2.32&2.61E-3&40.24&4.24E-3\\ 
C^3 &4096 &4.56&4.99E-3  &160.87&1.66E-2\\ \hline
\end{array} 
\]
\caption{The storage and time of matrix-vector multiplication for 1d nonlocal diffusion in heterogeneous medium. Kernel $\gamma(s)$ is $C^3$ away from origin. $\del_0=1/4$. The relative precision of HIF approximation is set to be $\epsilon=10^{-8}$.} \label{table:4}
\end{table}

\subsection{Two-dimensional test} \label{sec2:2d}
In this section, we perform 2d tests on the domain $\Om=(0,1)\times (0,1)$ to show the effectiveness of the HIF code. 
We choose $C(x,y)\equiv 1$ and $\del(x,y)\equiv\del$. 
The homogeneous medium has the same the horizon radius $\del$ for every point. 
The kernel $\gamma(s)$ is the same as before and has $O(1/|s|)$ singularity at $0$ and is assumed to be smooth expect at origin and $\pm1$.

In Table \ref{table:5}, the storage and the time of matrix-vector multiplication are listed for the original matrix $A$ as well as its compression after using the HIF method. 
The original matrix $A$ is again stored using the sparse matrix representation
and the MATLAB built-in matrix-vector multiplication is used to record the computation time in the right column of Table  \ref{table:5}.
The relative precision is chosen to be $1E-3$ in the 2d test. 
We remark that  in the case of $\del=1/2$, the critical improvement of the computational efficiency happens between
discontinuous kernel functions and $C^0$ kernel functions as a result of the choice of the precision parameter $\epsilon$.
Similarly as in 1d, higher critical regularity is observed if $\epsilon$ is getting smaller. 
 
  \begin{table}[htbp]
\[
 \begin{array}{|cc|cc|cc|}
\hline
 &  & \multicolumn{2}{|c|}{\text{HIF}} &   \multicolumn{2}{|c|}{\text{Original}}\\ \hline
\text{kernel} & \delta & \text{Memory (MB)}& T_{\text{multi}} \text{(sec)} &  \text{Memory (MB)}& T_{\text{multi}} \text{(sec)} \\ 
C^{-1} & 1 & 161.71 & 2.12E-2 & 248.24 & 2.50E-2\\
C^{0} & 1 & 33.01 & 8.93E-3 & 248.24 & 2.45E-2\\
C^{1} & 1 & 32.08 & 8.83E-3 & 248.24 & 2.37E-2\\
C^{2} & 1 & 31.50 & 8.59E-3 & 248.24 & 2.46E-2\\
C^{3}&1 &31.16&8.76E-3 & 248.24&2.42E-2\\ \hline
C^{-1} & 1/2& 264.63 &3.20E -2& 89.51 &9.41E-3\\
C^{0} & 1/2& 41.25 &1.04E -2& 89.51 &9.67E-2\\
C^{1} & 1/2& 30.98 &8.98E -3& 89.51 &1.00E-2\\
C^{2} & 1/2& 28.62 &8.86E -3& 89.51 &1.01E-2\\
C^{3}&1/2 &27.16&8.59E-3&89.51 &1.03E-2\\ \hline
C^{-1} & h &7.85&5.86E-3&0.49&8.17E-5\\
C^{0} & h &7.85&5.73E-3&0.49&8.05E-5\\
C^{1} & h &7.85&5.93E-3&0.49&8.78E-5\\
C^{2} & h &7.85&6.01E-3&0.49&8.68E-5\\
C^{3} & h &7.85&7.51E-3&0.49&9.40E-5\\ \hline
\end{array} 
\]

\caption{The storage and time of matrix-vector multiplication for 2d nonlocal diffusion in homogeneous  medium. The domain $\Om=(0,1)^2$ is discretized using $N=64\times 64$ points. The mesh size $h=1/64$. The relative precision of HIF approximation is set to be $\epsilon=10^{-3}$. } \label{table:5}
\end{table}

Table \ref{table:6} shows the storage and time for different value of $N$ with horizon parameter  fixed to be $1/2$. 
The discontinuous ($C^{-1}$) kernel and the $C^3$ kernel are used in the simulation. 
Here we remark that we choose the $C^{-1}$ kernel to be compared with the $C^3$ kernel because $C^{-1}$
is right before the critical improvement of efficiency happens with the increase of regularity as seen in Table \ref{table:5}. This is a result of the particular choice
of relative precision $\epsilon=1E-3$. 
As $\epsilon$ becomes smaller, higher regularity is needed in order for the critical improvement of efficiency to happen. 
We observe that the complexity of the matrix-vector multiplication for the original matrix $A$ without compression is approximately $O(N^2)$.
HIF also gives $O(N^2)$ computation time for matrix-vector multiplication with the discontinuous kernel being used, while 
in the case of $C^3$ kernel, HIF performs much better and the observed complexity for matrix-vector multiplication is about $O(N \log(N))$. 

 \begin{table}[htbp] 
\[
 \begin{array}{|cc|cc|cc|}
\hline
& &   \multicolumn{2}{|c|}{\text{HIF}} &   \multicolumn{2}{|c|}{\text{Original}}\\ \hline
\text{kernel} & N  & \text{Memory (MB)}&T_{\text{multi}} \text{(sec)}&  \text{Memory (MB)}& T_{\text{multi}} \text{(sec)}\\ 
C^{-1} &8^2   &0.07&1.18E-4&0.03&3.15E-5\\
C^{-1} &16^2 &1.05&2.29E-4&0.43&7.43E-5\\ 
C^{-1} &32^2&16.79&2.71E-3&6.01&8.19E-4\\ 
C^{-1} & 64^2 &264.28&3.29E-2&89.51&1.11E-2\\ \hline  
C^3& 8^2  &0.07&1.70E-4&0.03&3.24E-5\\
C^3& 16^2 &0.62&3.91E-4&0.43&7.49E-5\\ 
C^3& 32^2&4.19&1.73E-3&6.01&7.70E-4\\ 
C^3&  64^2 &27.21&9.05E-3  &89.51&1.09E-2\\ \hline
\end{array} 
\]
\caption{The storage and time of matrix-vector multiplication for 2d nonlocal diffusion in homogeneous medium. $\del=1/2$. The relative precision of HIF approximation is set to be $\epsilon=10^{-3}$.} \label{table:6}
\end{table}

To sum up, the simulations in Sections \ref{sec2:homo}, \ref{sec2:hetero} and \ref{sec2:2d} show that the popular choice of the kernel $\gamma$ in \eqref{kernelgamma} is problematic 
in the interest of fast computation since it has non-smooth truncation at $\pm1$. 
On the other hand, when the regularity of the kernel gets better the HIF fast algorithm becomes more effective.   
Those observations motivate us to propose the splitting method that deals kernels with non-smooth truncation in the simulations of the peridynamic model. 

\section{The splitting of singularities} \label{sec:3}
\begin{figure}[htbp] 
  \centering
\includegraphics[width=5.5cm]{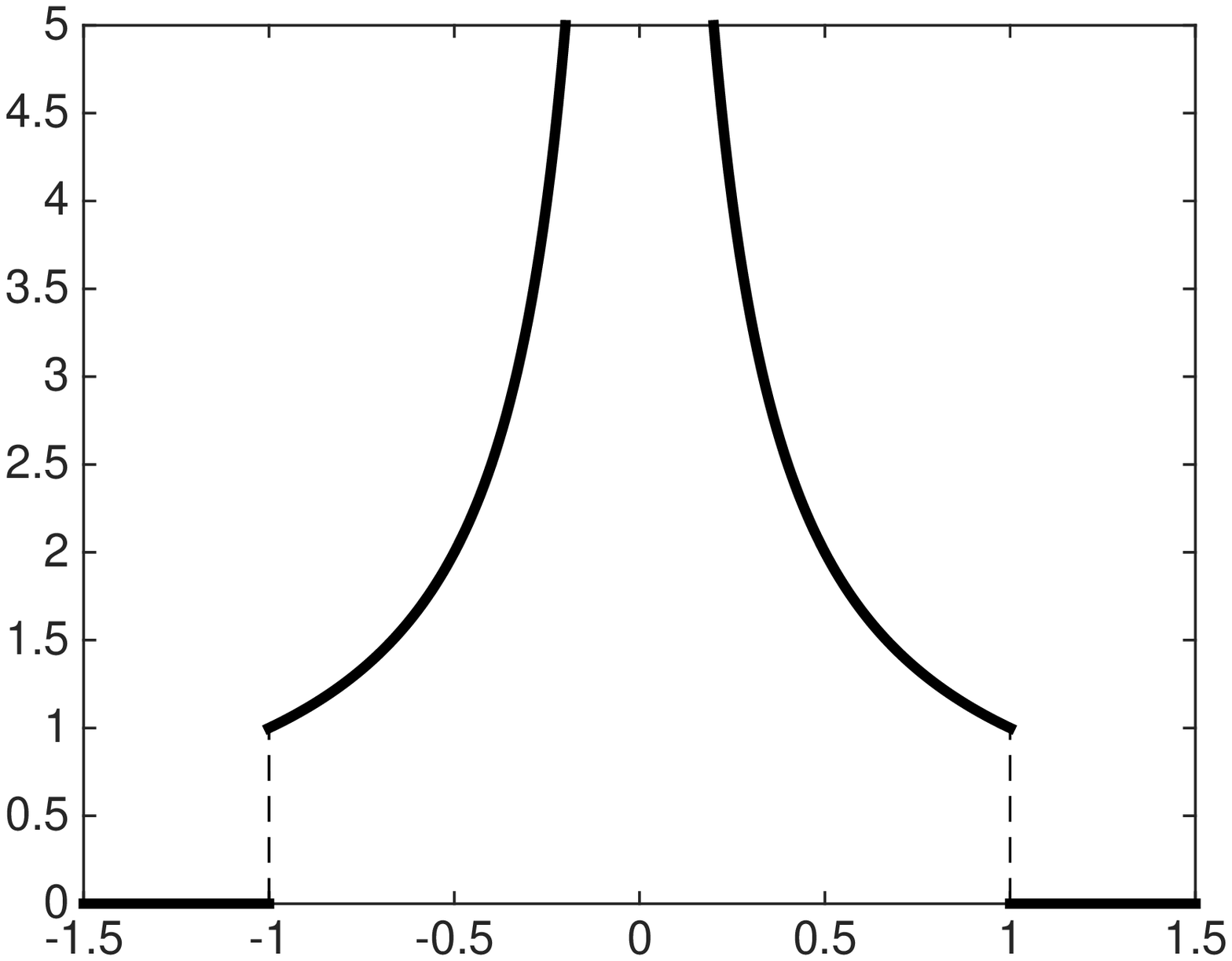}\\
\includegraphics[width=5.5cm]{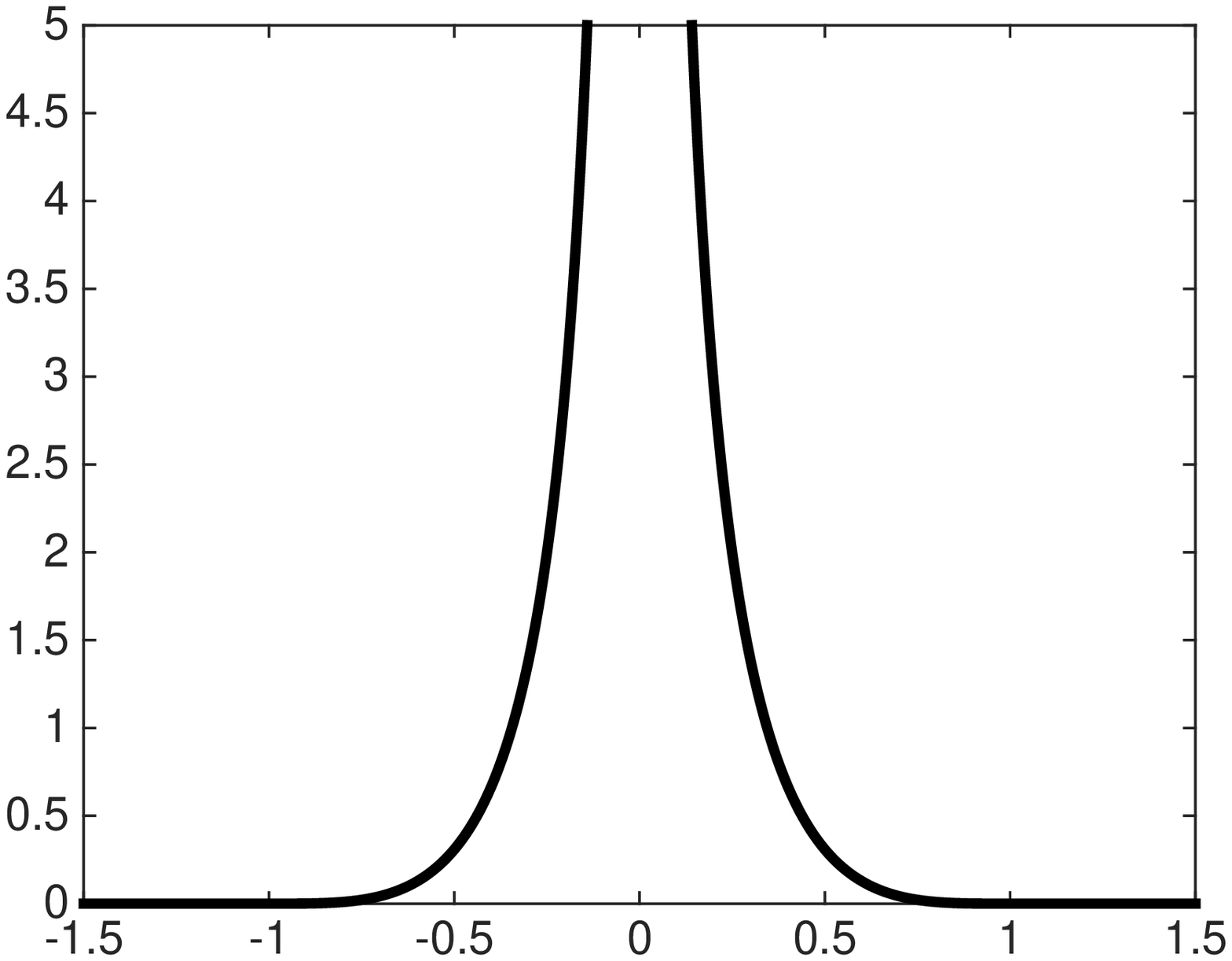}
\includegraphics[width=5.5cm]{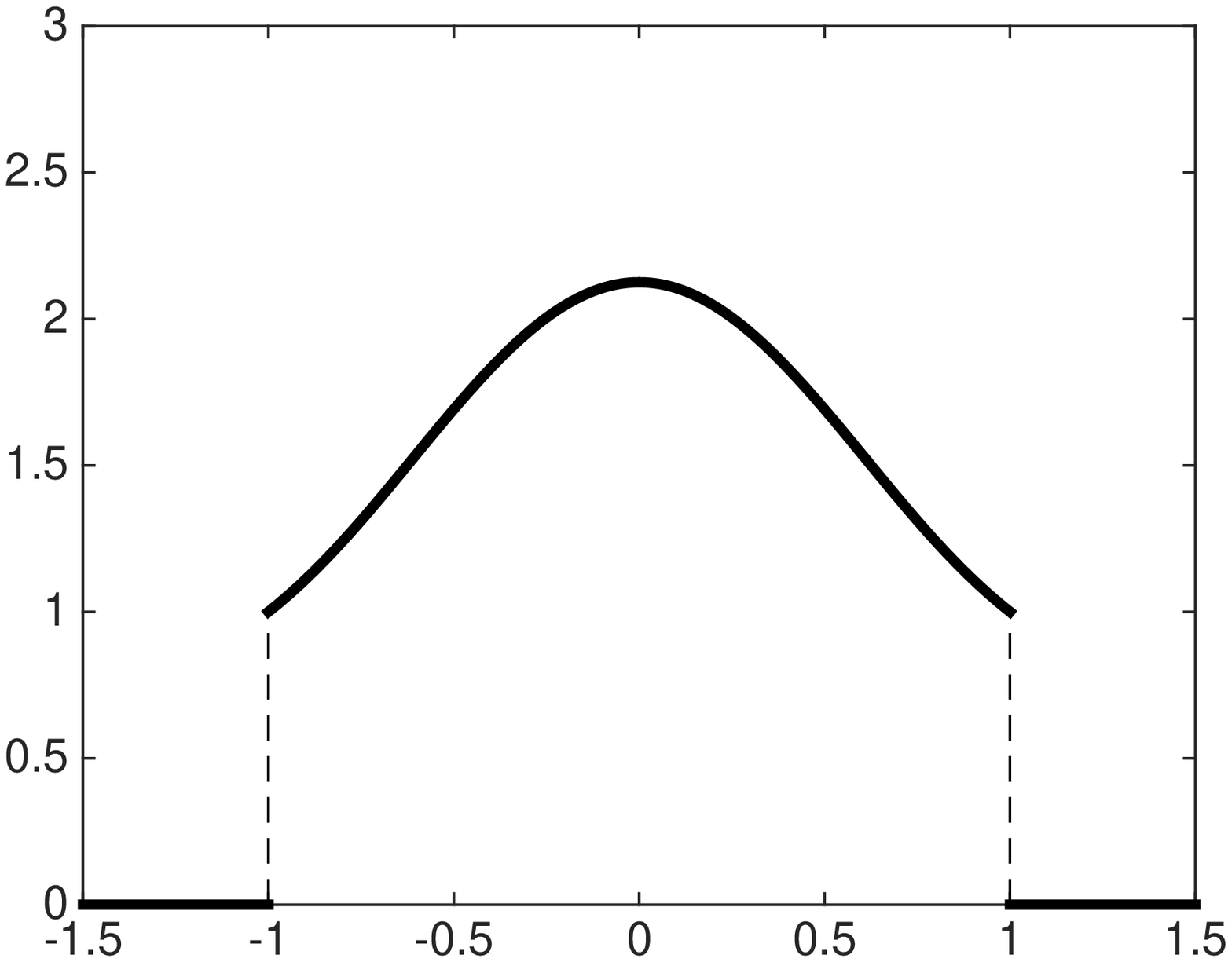}
\caption{The kernel $\gamma(s)$ (top) splits into $\kappa(s)$ (bottom left) and $p^3(s)\chi(|s|<1)$ (bottom right).} \label{fig:1}
    \end{figure}

Popular choices of kernels in the peridynamic simulations contain two types of singularities: the singularity at origin and 
the non-smooth truncation at the boundary  $\partial \cH_{\bx}$. The first type singularity is handled by the HIF algorithm through the idea of far-field compression. In 1d, the compression is performed through a block matrix with small block sizes near the diagonal. The algorithm is thus effective when the kernel is sufficiently smooth away from origin, as we have seen in the previous numerical examples.
Next, we propose to split the two types of singularities and treat them respectively. 
In Fig \ref{fig:1}, the left plot is the kernel $\gamma (|s|)= \frac{1}{|s|}\chi{(|s|<1)}$. 
The kernel splits into two parts, the middle plot which has singularity at zero but $C^3$ away from zero,
and the right plot which is discontinuous at $s=\pm1$ but smooth in the interior. 
In practice, the kernel $\gamma(|s|)$ is written into
\[
\gamma(|s|)= \kappa(|s|)+ p^{2K}(s)\chi(|s|<1) \,,
\] 
where $ \kappa(|s|)$ is $C^K$ regular away from zero, and $p^{2K}$ is an even polynomial of degree $2K$ that matches $\gamma$ and its derivatives at $|s|=1$ up to the $K$-th order. 
For the choice of $\gamma (|s|)= \frac{1}{|s|}\chi{(|s|<1)}$, $p^{2K}$ is explicitly given by the following formulas for $K=0,1,2,3$. 
\[
p^{2K}(s)=
\left\{
\begin{aligned}
&1 \quad &K=0 \\
&\frac{3}{2}-\frac{1}{2}s^2 \quad &K=1\\
&\frac{15}{8}-\frac{5}{4}s^2+\frac{3}{8}s^4 \quad& K=2\\
&\frac{17}{8}-2s^2+\frac{9}{8}s^4-\frac{1}{4}s^6 \quad & K=3
\end{aligned}
\right. 
\]

As suggested by Section \ref{sec:2}, the nonlocal operator corresponding to the kernel $ \kappa(s)$
can be treated by the HIF fast algorithm effectively. So we only need to deal with the second part
corresponding to the kernel $p^{2K}(s)\chi(|s|<1)$. 
Taking $\Om=[0,1]^d $ in the $d$-dimensional Euclidean space, we consider the 
following nonlocal diffusion problem with Dirichlet type boundary
\beq \label{nonlocal:polynomial}
\left\{
\begin{aligned}
&-\cL^P u (\bx) := -\int_{\cH_{\bx}} P(\bx,\by)(u(\by)-u(\bx)) d\by =f(\bx) &\quad x\in\Om\\
& u(\bx) = 0  &\quad x\in \Om^c
\end{aligned}
\right.
\eeq
where $\cH_{\bx}$ is some neighborhood of $\bx$, and $P$ is given by
\beq \label{kernel:polynomial}
P(\bx,\by)=\frac{C(\bx,\by)}{\del(\bx,\by)^{d+2}} p^{2K}(\frac{|\by-\bx|}{\del(\bx,\by)})\,.
\eeq
For convenience of illustration,  we take $C(\bx,\by)\equiv C(\bx)$, $\del(\bx,\by)\equiv \del(\bx)$, $K=0$ and $p^0\equiv 1$ in the next. 
The generalization to $p^{2K}$ for $K>0$ is straightforward. Taking into account of the boundary condition in \eqref{nonlocal:polynomial},
the nonlocal operator is to be evaluated at each $x$:  
\beq \label{nonlocalP}
\cL^P u (\bx)  = \frac{C(\bx)}{\del(\bx)^{d+2}} \left( \int_{\cH_{\bx}\cap \Om} u(\by)d\by -  |\cH_{\bx}| u(\bx) \right) \,.
\eeq

We therefore need only to find a fast algorithm compute $ \int_{\cH_{\bx}\cap \Om} u(\by)d\by$. 
Our algorithm starts by the hierarchical decomposition of the space $\Om$, which is similar to the one in standard 
FMM \cite{GrRo87}. 
Assume that the domain $\Om$ is uniformly  $N=n^d$ be the total number of discretization nodes. 
Let $L = \log_2(n)$, then from level $1$ to $L$, the computation domain $\Om=[0,1]^d$
is hierarchically subdivided into panels.
Each panel in the $l$-th level can be represented by one of the cubes $\prod_{j=1}^d [\frac{k_j}{2^l}, \frac{k_j+1}{2^l}]$ , $k_j=0, 1,\cdots, 2^l-1$.

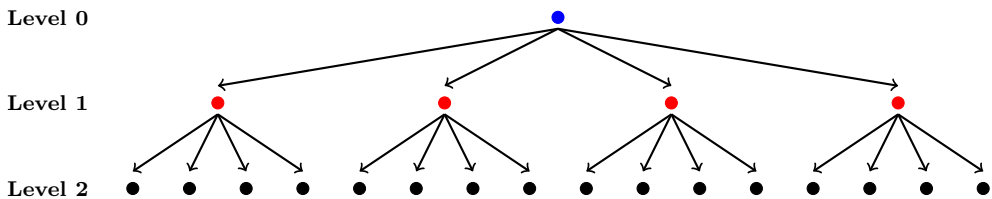
\begin{figure} [htbp]
    \centering
      \begin{tikzpicture}[scale=0.75]
      
      \node at (-9, 5) {\footnotesize \textbf{Level 0}};
       \node at (-9, 3. 5) {\footnotesize \textbf{Level 1}};
        \node at (-9, 2) {\footnotesize \textbf{Level 2}};

       \draw(0,5) circle(3pt) [blue, fill=blue];
       \draw [->, thick] (0,4.8)--(-6,3.8); 
        \draw [->, thick] (0,4.8)--(-2,3.8); 
         \draw [->, thick] (0,4.8)--(2,3.8); 
          \draw [->, thick] (0,4.8)--(6,3.8); 
          
      \draw(-6,3.5) circle(3pt) [red, fill=red];
      \draw(-2,3.5) circle(3pt) [red, fill=red];
      \draw(2,3.5) circle(3pt) [red, fill=red];
      \draw(6,3.5) circle(3pt) [red, fill=red];
      
      \draw [->, thick] (-6,3.3)--(-7.5,2.3); 
       \draw [->, thick] (-6,3.3)--(-6.5,2.3); 
        \draw [->, thick] (-6,3.3)--(-5.5,2.3); 
         \draw [->, thick] (-6,3.3)--(-4.5,2.3); 
         
            \draw(-7.5,2) circle(3pt) [black, fill=black];
      \draw(-6.5,2) circle(3pt) [black, fill=black];
      \draw(-5.5,2) circle(3pt) [black, fill=black];
      \draw(-4.5,2) circle(3pt) [black, fill=black];
         
         \draw [->, thick] (-2,3.3)--(-3.5,2.3); 
       \draw [->, thick] (-2,3.3)--(-2.5,2.3); 
        \draw [->, thick] (-2,3.3)--(-1.5,2.3); 
         \draw [->, thick] (-2,3.3)--(-0.5,2.3); 
         
               \draw(-3.5,2) circle(3pt) [black, fill=black];
      \draw(-2.5,2) circle(3pt) [black, fill=black];
      \draw(-1.5,2) circle(3pt) [black, fill=black];
      \draw(-0.5,2) circle(3pt) [black, fill=black];
         
          \draw [->, thick] (2,3.3)--(0.5,2.3); 
       \draw [->, thick] (2,3.3)--(1.5,2.3); 
        \draw [->, thick] (2,3.3)--(2.5,2.3); 
         \draw [->, thick] (2,3.3)--(3.5,2.3); 
                  
      \draw(0.5,2) circle(3pt) [black, fill=black];
      \draw(1.5,2) circle(3pt) [black, fill=black];
      \draw(2.5,2) circle(3pt) [black, fill=black];
      \draw(3.5,2) circle(3pt) [black, fill=black];
         
       \draw [->, thick] (6,3.3)--(4.5,2.3); 
       \draw [->, thick] (6,3.3)--(5.5,2.3); 
        \draw [->, thick] (6,3.3)--(6.5,2.3); 
         \draw [->, thick] (6,3.3)--(7.5,2.3); 
                           
      \draw(4.5,2) circle(3pt) [black, fill=black];
      \draw(5.5,2) circle(3pt) [black, fill=black];
      \draw(6.5,2) circle(3pt) [black, fill=black];
      \draw(7.5,2) circle(3pt) [black, fill=black];
    \end{tikzpicture}
          \caption{Tree structure induced by the hierarchical subdivision of $[0,1]^2$.} \label{fig:tree}
          \end{figure}

    The hierarchical subdivision of $\Om$ forms a tree structure. Fig \ref{fig:tree} is the quadtree formed by the subdivision of 
    the two-dimensional domain $[0,1]^2$. The blue dot in Fig \ref{fig:tree} represents the original domain $[0,1]^2$,
    and it has four children $\{ [\frac{i}{2}, \frac{i+1}{2}]\times[\frac{j}{2}, \frac{j+1}{2}] \}_{i,j=0}^1$ represented by the red dots in level $1$. 
    The subdivision of the $4$ domains in level $1$ results in $16$ subdomains in level $2$ represented by the black dots. 
    The hierarchical subdivision is performed until reaching the $L$-th level. 
     Similarly, in 1d a binary tree will be formed by the subdivision of $[0,1]$ and  in 3d an octree will be formed by the subdivision of $[0,1]^3$. 
    
We now present the algorithm. It contains the following $3$ steps that will be explained in detail. 

\begin{enumerate} 
{\itshape 
\item Initialization. Traverse the tree from top to bottom to obtain the decomposition of $\cH_{\bx}\cap \Om$ for each grid point $\bx$. 
\item Traverse the tree from bottom to top to obtain the partial sums $\sum_{\bx_k \in Q_l^j} u(\bx_k)$ for each panel $ Q_l^j$ in the $l$-th level.
\item Approximate the integral  $\int_{\cH_{\bx}\cap \Om} u(\by)d\by$ by using (1) and (2). \\
} 
\end{enumerate}

\textbf{Step 1}. {\itshape Initialization step}. 
For each $\bx$ on the grid point, decompose the domain of integration $\cH_{\bx}\cap \Om$
 using the panels in $L_l$ ($l=0,1,\cdots, L$). The decomposition of $\cH_{\bx}\cap \Om$ for each $\bx$ is stored 
 for reuse. More specifically, we traverse the tree from top to bottom to determine whether a panel $Q_l^j$ from level $l$ is included or not,
and then express $\cH_{\bx}\cap \Om = \bigcup_{l=0}^L\bigcup_{j\in \mathcal{I}_l(\bx) } Q_l^j $, where $\mathcal{I}_l(\bx)$ is the set of indices at the $l$-th to be included in the decomposition. 
 We may do this recursively by calling the function recur($\Om$), where recur is defined in Algorithm \ref{alg:recur} below. 

  \begin{algorithm}[htbp]

function recur($Q$)\\
\uIf{$Q$ is fully contained in  $ \cH_{\bx}\cap \Om$}{
    $Q$ is included\;
     \textbf{return} 1; 
  }
  \uElseIf{$Q$ intersects with $(\cH_{\bx}\cap \Om)$  nontrivially}{
    \eIf{ $Q$ is not a leaf}{
    run recur($Q[k]$) from $k=1$ to $k=2^d$, where $Q[k]$ is the $k$th children of $Q$\;
    }
    { $Q$ is included\;
     \textbf{return} 1;  }
  }
  \Else{
     \textbf{return} 0\;
  }
    \caption{The recursive function.}  \label{alg:recur}
 \end{algorithm}

 Let us now discuss on the complexity of performing the recur function when $\cH_{\bx} $ is a $d$-dimensional ball centered at $\bx$. The main cost of running the recur function comes from the intersection test of balls and cubes. To test whether a cube $Q$ is fully contained in a ball, we only need to test whether the total of $2^d$ corners of $Q$ are all contained in the ball. This requires a total of  $O(2^d)$ operations, where the constant in $O(2^d)$ is independent of the dimension $d$. Now to test whether a cube $Q$ has nontrivial intersection with a ball, we use the 
 algorithm given in \cite{Arvo2013}, which is an $O(d)$ algorithm that determines whether an axis-aligned bounding box (AABB) intersects
 a ball. See the function defined in Algorithm \ref{alg:intersectiontest}, where the inputs are three vectors $\textbf{\text{Bmin}},\textbf{\text{Bmax}}, \textbf{\text{C}}$ and a positive number $r$. The vector $\textbf{\text{Bmin}}\in\R^d$ stores the minima of the AABB for each axis, $\textbf{\text{Bmax}}\in \R^d$ stores the maxima of the AABB for each axis, and $\textbf{\text{C}}\in\R^d$ and $r$ are the center and radius of the ball resepctively. Combining the above discussions, we know that each time calling the recur function
 requires $C \cdot 2^d$ operations, where $C$ is independent of dimension.   
 
 \begin{algorithm}[htbp]
function doesCubeIntersectBall (\textbf{\text{Bmin}}, \textbf{\text{Bmax}}, \textbf{\text{C}}, r)\\
$dmin$ = 0; \\
\For{$i = 1:d$}{
     \uIf{ \rm $\textbf{\text{Bmin}} [i]>\textbf{\text{C}} [i]$}{
     $dmin =dmin+ (\textbf{\text{Bmin}} [i]- \textbf{\text{C}} [i])^2$; 
     }
   \uElseIf{\rm $\textbf{\text{Bmax}} [i]<\textbf{\text{C}} [i]$}{
     $dmin =dmin+ (\textbf{\text{Bmax}} [i]- \textbf{\text{C}} [i])^2$; 
     }
  }  
   \uIf{$dmin<r^2$}{
     \textbf{return} TRUE ;
     } 
      \Else{
    \textbf{return} FALSE ;
  }
  
    \caption{The intersection test of $d$-dimensional cubes and balls.}  \label{alg:intersectiontest}
 \end{algorithm}

 \begin{figure}[htbp]
\centering
\includegraphics[width=8cm]{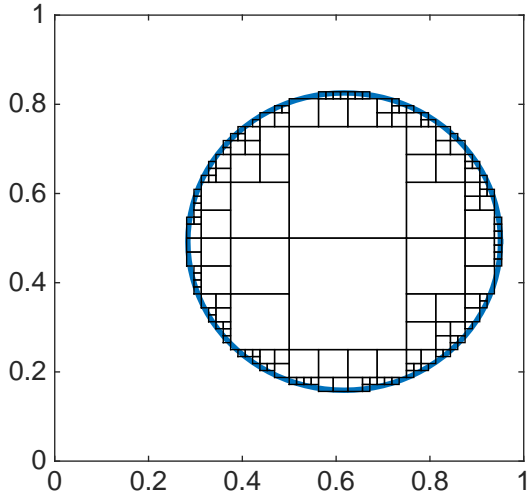}
\vspace{-0.5cm}
\caption{The decomposition of a circular region into panels in $\{ L_i \}$}
 \label{fig:decomposition}
\end{figure}

Fig \ref{fig:decomposition} shows an example of decomposing a circular region in 2d. 
Now the total complexity of Step 1 is obviously equal to $C\cdot 2^d $ times the total number of times the recur function is called.
Now for each discretization node $\bx_i$ ($i=1, 2, \cdots, N$), let $\mathscr{N}_d(\bx_i)$
denote the number of times the recur function is called in order to decompose the domain $\cH_{\bx_i}\cap \Om$.
The total number of times the recur function is called is then given by $\sum_{i=1}^N \mathscr{N}_d(\bx_i)$.  
Now observe that $\mathscr{N}_d (\bx_i)$ is also equal to the number of panels that intersect with the boundary $\partial \cH_{\bx_i}$ non-trivially, namely
\[
\mathscr{N}_d (\bx_i) = \# \bigcup_{l=0}^{L}\{ Q_l^j \in L_l : (Q_l^j)^{\mathrm{o}} \cap \partial \cH_{\bx_i} \neq \emptyset \} =\sum_{l=0}^L \mathscr{N}^l_d (\bx_i)\,,
\]
where  $\mathscr{N}^l_d (\bx_i)$ denotes the number of panels  in the $l$-th level that intersect with $\partial \cH_{\bx_i}$ non-trivially. 
Assume that  for each point $\bx_i$, the set $\partial \cH_{\bx_i}$ has uniform bounded mean curvature, 
then $\mathscr{N}^l_d (\bx_i)$  can be estimated by the following formula 
\[
\mathscr{N}^l_d (\bx_i)= 
\left\{ 
\begin{aligned}
&O(1) \quad &\text{for } d=1 \,;\\
&O\left(\frac{N}{(2^d)^l}\right)^{1 -1/d}  \quad &\text{for } d\geq2 \,.
\end{aligned}
\right.
\]
To sum up, the total complexity of Step 1, which equals $C\cdot 2^d\sum_{i=1}^N \sum_{l=0}^L \mathscr{N}^l_d (\bx_i) $ is given by 
\beq \label{eq:complexity}
\left\{ 
\begin{aligned}
O(N\log N) \quad \text{for } d=1 \,;\\
C_d \cdot N^{2 -1/d}   \quad \text{for } d\geq2 \,,
\end{aligned}
\right.
\eeq
where the dimensional dependent constant $C_d$ is of order $O(2^d)$.

\textbf{Step 2}. {\itshape Compute the the partial sums for each tree node} . For a given vector $\{u(\bx_i)\}_{i=1}^N$, we assign the value $u(\bx_i)$ to a leaf node.   
Then traverse the tree bottom to top, we assign each parent node the value of the sum of all its children. 
The completion of this step gives us the partial sums $\sum_{\bx_k \in Q_l^j} u(\bx_k)$ for every panel $Q_l^j$. 
The partial sums obtained  can then be used to approximate the integral $\int_{Q_l^j} u(\by)d\by$. We remark that if the kernel function given by \eqref{kernel:polynomial} is generated from the polynomial $p^{2K}$ with $K>0$, we then need to compute the partial sums in the form of
$\sum_{x_k \in Q_l^j} x_k^m u(x_k)$ for all $m=0,1,\cdots 2K$, so that they can be used to approximate the integral $\int_{Q_l^j} P(\bx, \by) u(\by)d\by$.

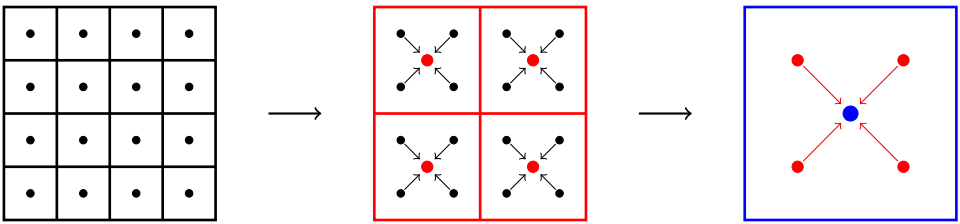
\begin{figure}[htbp] 
    \centering
      \begin{tikzpicture}[scale=0.7]
      \draw [line width=1, blue ] (-2+14, -2) rectangle (2+14, 2);
       \draw(14,0) circle(4pt) [blue, fill=blue];
                \draw(6+7,-1) circle(3pt) [red, fill=red];
             \draw[->,red, shorten <=3pt, shorten >=5pt] (6+7,-1)  -- (14,0) ; 
           \draw(6+7,1) circle(3pt) [red, fill=red];
                        \draw[->,red, shorten <=3pt, shorten >=5pt] (6+7,1)  -- (14,0) ; 
             \draw(8+7,-1) circle(3pt) [red, fill=red];
                          \draw[->,red, shorten <=3pt, shorten >=5pt] (8+7,-1)  -- (14,0) ; 
               \draw(8+7,1) circle(3pt) [red, fill=red];
                            \draw[->,red, shorten <=3pt, shorten >=5pt] (8+7,1)  -- (14,0) ; 
      \draw [->, thick] (3,0)--(4,0);
      
      \draw [line width=1, red ] (5, -2) rectangle (9, 2);
      \draw [line width=1, red] (5,0) -- (9,0);
      \draw [line width=1, red] (7,-2) -- (7,2);
         \draw(6,-1) circle(3pt) [red, fill=red];
           \draw(6,1) circle(3pt) [red, fill=red];
             \draw(8,-1) circle(3pt) [red, fill=red];
               \draw(8,1) circle(3pt) [red, fill=red];
                        \draw(12.5-7,-0.5) circle(2pt) [black, fill=black];
                          \draw[->,shorten <=2pt, shorten >=4pt] (12.5-7,-0.5)  -- (6,-1) ; 
                         \draw(12.5-7,0.5) circle(2pt) [black, fill=black];
         \draw[->,shorten <=2pt, shorten >=4pt] (12.5-7,0.5)  -- (6,1) ; 
           \draw(12.5-7,-1.5) circle(2pt) [black, fill=black];
                  \draw[->,shorten <=2pt, shorten >=4pt] (12.5-7,-1.5)  -- (6,-1) ; 
          \draw(12.5-7,1.5) circle(2pt) [black, fill=black];
            \draw[->,shorten <=2pt, shorten >=4pt] (12.5-7,1.5)  -- (6,1) ; 
             \draw(13.5-7,-0.5) circle(2pt) [black, fill=black];
                           \draw[->,shorten <=2pt, shorten >=4pt] (13.5-7,-0.5)  -- (6,-1) ; 
          \draw(13.5-7,0.5) circle(2pt) [black, fill=black];
            \draw[->,shorten <=2pt, shorten >=4pt] (13.5-7,0.5)  -- (6,1) ; 
           \draw(13.5-7,-1.5) circle(2pt) [black, fill=black];
              \draw[->,shorten <=2pt, shorten >=4pt] (13.5-7,-1.5)  -- (6,-1) ; 
          \draw(13.5-7,1.5) circle(2pt) [black, fill=black];
            \draw[->,shorten <=2pt, shorten >=4pt] (13.5-7,1.5)  -- (6,1) ; 
             \draw(14.5-7,-0.5) circle(2pt) [black, fill=black];
                           \draw[->,shorten <=2pt, shorten >=4pt] (14.5-7,-0.5)  -- (8,-1) ; 
          \draw(14.5-7,0.5) circle(2pt) [black, fill=black];
               \draw[->,shorten <=2pt, shorten >=4pt] (14.5-7,0.5)  -- (8,1) ; 
           \draw(14.5-7,-1.5) circle(2pt) [black, fill=black];
                \draw[->,shorten <=2pt, shorten >=4pt] (14.5-7,-1.5)  -- (8,-1) ; 
          \draw(14.5-7,1.5) circle(2pt) [black, fill=black];
                         \draw[->,shorten <=2pt, shorten >=4pt] (14.5-7,1.5)  -- (8,1) ; 
             \draw(15.5-7,-0.5) circle(2pt) [black, fill=black];
                  \draw[->,shorten <=2pt, shorten >=4pt] (15.5-7,-0.5)  -- (8,-1) ; 
          \draw(15.5-7,0.5) circle(2pt) [black, fill=black];
                         \draw[->,shorten <=2pt, shorten >=4pt] (15.5-7,0.5)  -- (8,1) ; 
           \draw(15.5-7,-1.5) circle(2pt) [black, fill=black];
                \draw[->,shorten <=2pt, shorten >=4pt] (15.5-7,-1.5)  -- (8,-1) ; 
          \draw(15.5-7,1.5) circle(2pt) [black, fill=black];
                         \draw[->,shorten <=2pt, shorten >=4pt] (15.5-7,1.5)  -- (8,1) ; 
       \draw [->, thick] (10,0)--(11,0);
       
         \draw [line width=1, black ] (12-14, -2) rectangle (16-14, 2);
         \draw [line width=1, black ] (12-14, 0)-- (16-14, 0);
         \draw [line width=1, black ] (12-14, 1)-- (16-14, 1);
         \draw [line width=1, black ] (12-14, -1)-- (16-14, -1);
           \draw [line width=1, black ] (13-14, -2)-- (13-14, 2);
             \draw [line width=1, black ] (14-14, -2)-- (14-14, 2);
               \draw [line width=1, black ] (15-14, -2)-- (15-14, 2);
         \draw(12.5-14,-0.5) circle(2pt) [black, fill=black];
          \draw(12.5-14,0.5) circle(2pt) [black, fill=black];
           \draw(12.5-14,-1.5) circle(2pt) [black, fill=black];
          \draw(12.5-14,1.5) circle(2pt) [black, fill=black];
             \draw(13.5-14,-0.5) circle(2pt) [black, fill=black];
          \draw(13.5-14,0.5) circle(2pt) [black, fill=black];
           \draw(13.5-14,-1.5) circle(2pt) [black, fill=black];
          \draw(13.5-14,1.5) circle(2pt) [black, fill=black];
             \draw(14.5-14,-0.5) circle(2pt) [black, fill=black];
          \draw(14.5-14,0.5) circle(2pt) [black, fill=black];
           \draw(14.5-14,-1.5) circle(2pt) [black, fill=black];
          \draw(14.5-14,1.5) circle(2pt) [black, fill=black];
             \draw(15.5-14,-0.5) circle(2pt) [black, fill=black];
          \draw(15.5-14,0.5) circle(2pt) [black, fill=black];
           \draw(15.5-14,-1.5) circle(2pt) [black, fill=black];
          \draw(15.5-14,1.5) circle(2pt) [black, fill=black];

          \end{tikzpicture}
          \caption{A 2d example of the process of computing partial sums.} \label{fig:subdivision}
 
         \end{figure}

Fig \ref{fig:subdivision} exemplifies the way partial sums are computed for a two-dimensional problem. 
Each of the black dots on the left figure represents a leaf node that contains the partial sum of the function in the little black box that 
includes it. 
Then the partial sums are passed to the coarser level where each parent node presented by the red dots stores the sum of the values from its four children. The total sum of the function on the whole domain is stored in the top tree node presented by the blue dot in the right figure in Fig \ref{fig:subdivision}.  Therefore, the computational cost of Step 2 depends only on the number of discretization nodes $N$ and the degree $2K$ of the polynomial $p^{2K}$. To sum up, the  complexity of Step 2 is given by $O(K N)$, where the constant in $O(K N)$ is independent of dimension.

\textbf{Step 3}.  For each discretization node $\bx_i$ ($i=1,2\cdots N$), use the decomposition of $\cH_{\bx_i}\cap \Om$ obtained in Step 1
 and the partial sums obtained in step 2 to approximate the integral
 \beq\label{eq:sum}
 \int_{\cH_{\bx_i}\cap \Om} u(\by)d\by = \sum_{l=0}^L \sum_{j\in \mathcal{I}_l(\bx_i) } \int_{Q_l^j} u(\by)d\by  \,.
 \eeq
  Finally,  the value of $\cL_\del^Pu(\bx_i)$ is obtained through the relation \eqref{nonlocalP}. We again remark that in the case of $K>0$ for the polynomial $p^{2K}$, 
 all the partial sums of the form $\sum_{x_k \in Q_l^j} x_k^m u(x_k) (m=0,1,\cdots, 2K)$ will be
 needed to approximate the desired integral. 
With the precomputed sums for the approximate integrals of the form $\int_{Q_l^j} u(\by)d\by $ in Step 2, the number of summations we need to take in \eqref{eq:sum} is essentially equal to
 $\mathscr{N}_d (\bx_i)$ defined in Step 1. By the calculations in Step 1 and taking into the account of the general case that $p^{2K}$ being a polynomial of degree $2K>0$,  complexity of Step 3 is then given by the following formulas,

\beq \label{eq:complexity3}
\left\{ 
\begin{aligned}
O(K N\log N) \quad \text{for } d=1 \,;\\
O( K N^{2 -1/d})   \quad \text{for } d\geq2 \,,
\end{aligned}
\right.
\eeq

and the constants in the above estimates are independent of dimension. 

\subsection{Numerical tests} We perform numerical tests based on the proposed algorithm.\\

 \noindent{\itshape Example 1}. In this example, 
we take $\Om=[0,1]^2$ in the 2-dimensional space. Assume $k=0$, $p^0(x)\equiv 1$, $\del(\bx)=\frac{1}{2}$ and $C(\bx)=1$. 
The nonlocal operator to be evaluated is given by
\[
\cL_\del^P u(\bx) = \frac{1}{\del^{4}}\int_{B_\del(\bx)} (u(\by)-u(\bx)) d\by = \frac{1}{\del^{4}}\int_{B_\del(\bx)\cap \Om} u(\by) d\by -\frac{|B_1|}{\del^2}u(\bx)\,.
\]

\begin{figure}[htbp]
\centering
\includegraphics[width=7cm]{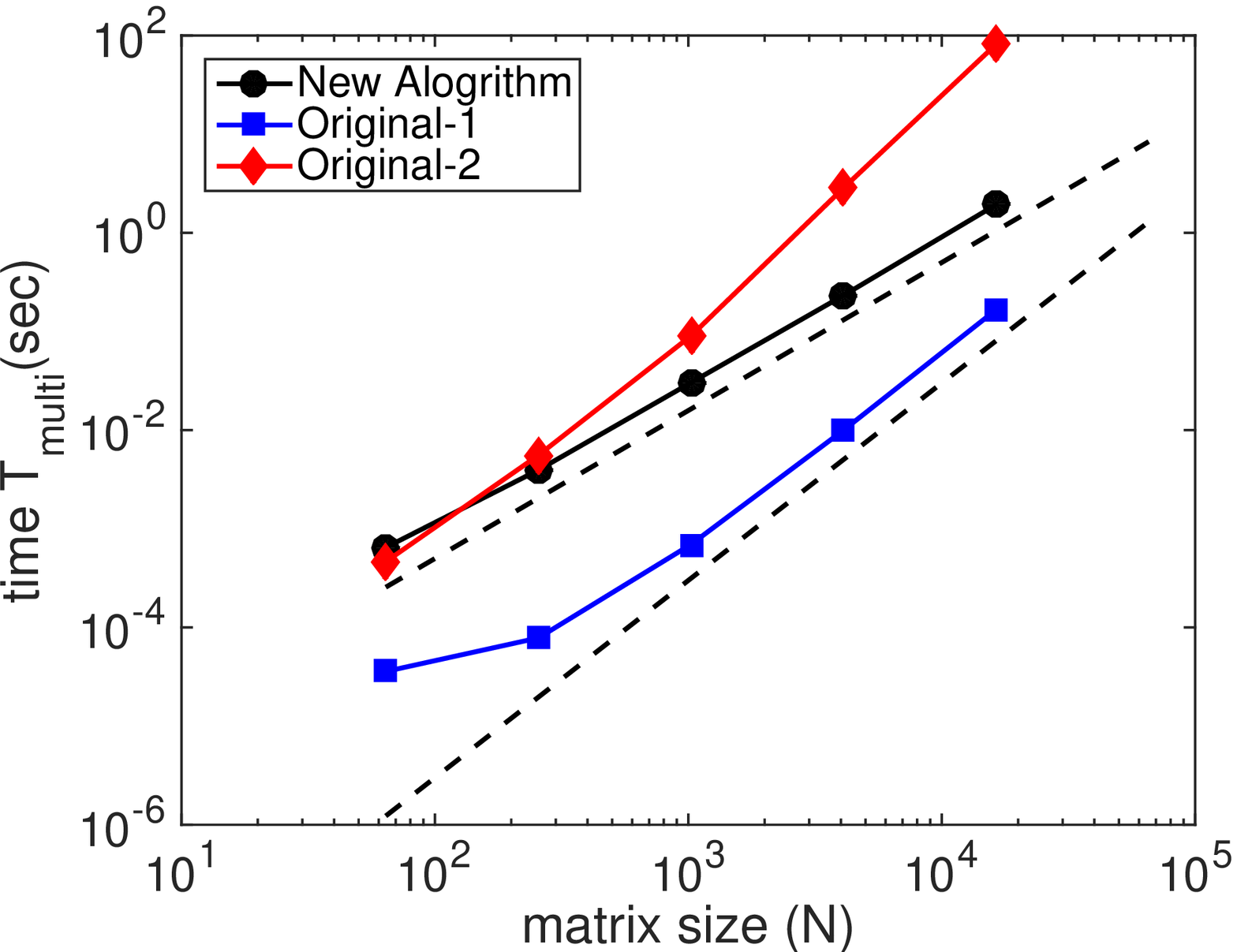}
\caption{Scaling results for Example 1, comparing the new algorithm with the two versions of original matrix-vector multiplication.
The reference scalings (dashed lines) of $O(N^2)$ and $O(N^{1.5})$ are also included.}
 \label{fig:scale1}
\end{figure}

Table \ref{table:sm2d} contains the experimental data of the new algorithm in comparison with two versions of matrix-vector multiplication
for the original matrix. The scaling results are shown in Fig.~\ref{fig:scale1}.
Since the initialization step (Step 1) is only needed to be performed once in an iterative method or time-dependent problem,
 the computation time recorded for the new algorithm is given by the cost for Step 2 and Step 3 combined together.  
The computation time of the new algorithm grows at the rate $O(N^{1.5})$, which verifies that the the formula given in \eqref{eq:complexity} for $d=2$. In comparison, the two versions of matrix-vector multiplication
for the original matrix shows $O(N^2)$ complexity. The first version (Original - 1) is performed with the MATLAB built-in matrix-vector multiplication. Because of the nature of MATLAB, the second version (Original -2),
which goes through the evaluation of the nonlocal operator at each  node $\bx_i$ $(i=1, 2,\cdots N)$ by using a for-loop, 
is also listed for a fair comparison 
with the new algorithm. 
 \\ 

  \begin{table}[htbp]
\[
 \begin{array}{|c|cc|cc|cc|}
\hline
&    \multicolumn{2}{|c|}{\text{New Algorithm}} &   \multicolumn{2}{|c|}{\text{Original - 1}}& \multicolumn{2}{|c|}{\text{Original - 2}} \\ \hline
N  &\text{ }T_{\text{multi}} \text{(sec) }&\text{ rate } &\text{ }T_{\text{multi}} \text{(sec) }&\text{ rate }&\text{ }T_{\text{multi}} \text{(sec) }&\text{ rate }\\ 
8^2 &6.47E-4& -- &3.58E-5&--&4.64E-4&--\\
16^2 &3.94E-3&1.30&7.83E-5&0.56&5.56E-3&1.79\\
32^3 &3.03E-2&1.47&6.94E-4&1.57&9.02E-2&2.01\\ 
64^2 &2.32E-1&1.47&9.99E-3&1.92&2.85E+0&2.50 \\
128^2 & 1.93E+0& 1.53&1.66E-1&2.03&8.14E+1&2.42\\ \hline
\end{array} 
\]
\caption{Computation time for evaluating the nonlocal operator in $2d$. $k=0$ and $p^0=1$.
 $T_{\text{multi}}$ under the new algorithm is the time for Step 2 and Step 3 combined.
 Original -1  uses the MATLAB built-in sparse matrix-vector multiplication. 
 Original -2 uses a for-loop to go through $\bx_i (i=1,2,\cdots,N)$ together with the MATLAB built-in vector-vector multiplication applied in each loop. } \label{table:sm2d}
\end{table}

\noindent{\itshape Example 2}.  To show that the algorithm can be applied to more generalized cases,
we perform the numerical test on the 1d domain $\Om=[0,1]$ and the kernel function $P(x, y)$ is given by \eqref{kernel:polynomial}
 with $C(x,y)=1$, $\del(x,y)=\del(x)$ as shown in Fig \ref{fig:delta}. We take $K=3$ and $p^{6}$ a polynomial of degree $6$. 
Now in Step 2 of the algorithm, we not only need to compute the partial sums $\sum_{x_k \in Q_l^j} u(x_k)$, but also
$\sum_{x_k \in Q_l^j} x_k^m u(x_k)$ for all $m=1,2,\cdots 6$ so that they can be used to approximate $\int_{\cH_{x}\cap \Om} P(x, y) u(y)dy$. 

\begin{figure}[htbp]
\centering
\includegraphics[width=7cm]{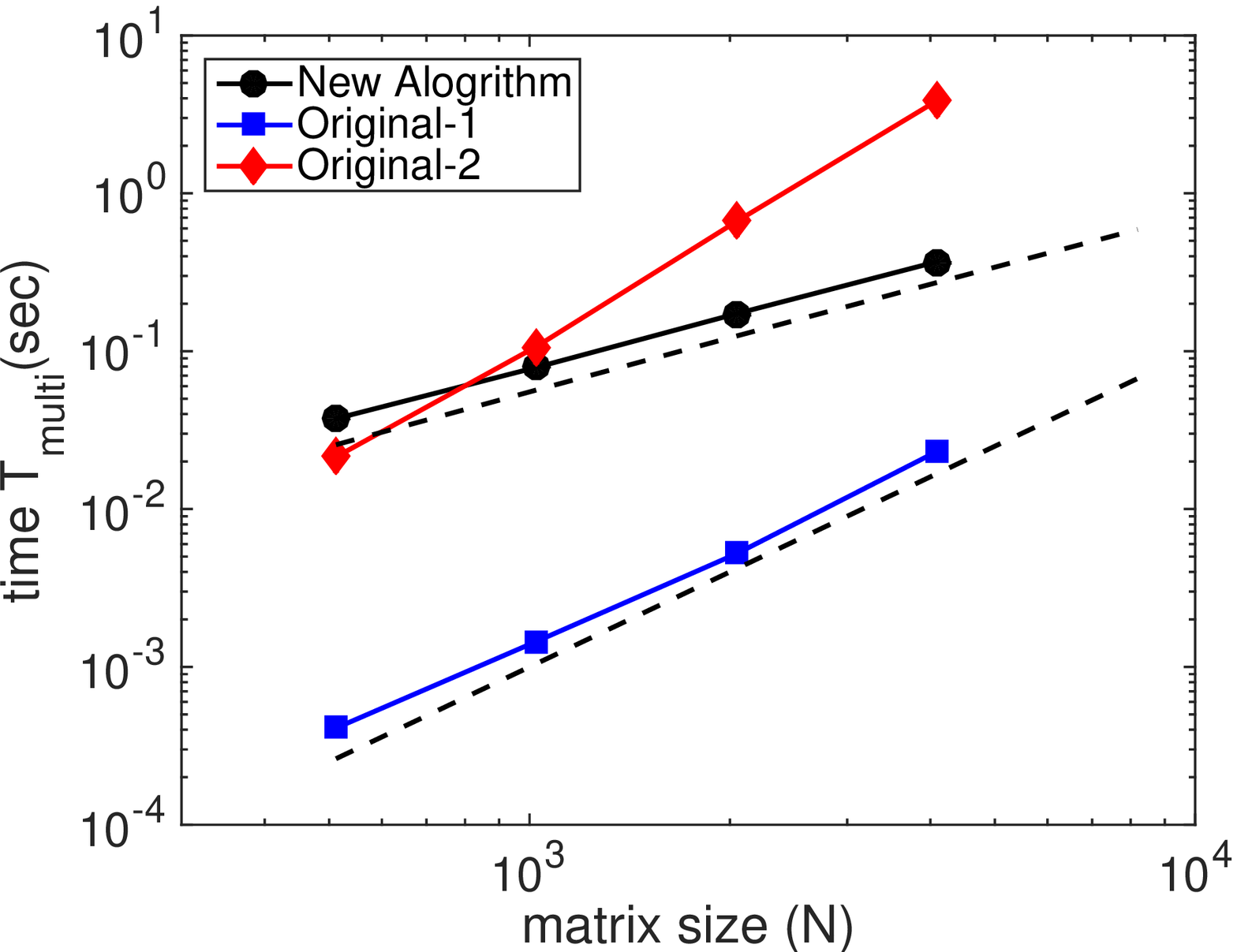}
\caption{Scaling results for Example 2, comparing the new algorithm with the two versions of original matrix-vector multiplication.
The reference scalings (dashed lines) of $O(N^2)$ and $O(N\log(N))$ are also included. }
 \label{fig:scale2}
\end{figure}

The computation time for the new algorithm is given in Table \ref{table:sm1d} in comparison with 
two versions of matrix-vector multiplication
for the original matrix, with scaling results shown in Fig.~\ref{fig:scale2}.
The computation time of the new algorithm grows at the rate which is slightly higher than $O(N)$, which relates well to our theoretical estimate of $O(N\log(N))$ for
$d=1$. 
\begin{table}[htbp]
\[
 \begin{array}{|c|cc|cc|cc|}
\hline
&    \multicolumn{2}{|c|}{\text{New Algorithm}} &   \multicolumn{2}{|c|}{\text{Original - 1}}& \multicolumn{2}{|c|}{\text{Original - 2}} \\ \hline
N  &\text{ }T_{\text{multi}} \text{(sec) }&\text{ rate } &\text{ }T_{\text{multi}} \text{(sec) }&\text{ rate }&\text{ }T_{\text{multi}} \text{(sec) }&\text{ rate }\\ 
512  &3.74E-2 & --&4.10E-4&--&2.14E-2&--\\
1024 &7.90E-2& 1.08&1.45E-3&1.82&1.07E-1&2.32\\ 
2048&1.73E-1&1.13&5.22E-3&1.85&6.65E-1&2.64\\ 
4096 &3.69E-1 &1.10&2.32E-2&2.15&3.87E+0&2.54\\ \hline
\end{array} 
\]

\caption{Computation time for evaluating the nonlocal operator in 1d. $k=3$ and $p^3$ is an even polynomial of degree $6$. $\del(x)$ is given
by Fig \ref{fig:delta} with $\del_0=1/4$.  $T_{\text{multi}}$ under the new algorithm is the time for Step 2 and Step 3 combined.
 Original -1  uses the MATLAB built-in sparse matrix-vector multiplication. 
 Original -2 uses a for-loop to go through $\bx_i (i=1,2,\cdots,N)$ together with the MATLAB built-in vector-vector multiplication applied in each loop.} \label{table:sm1d}
\end{table}

\section{Conclusion}
In this work, we develop and investigate fast algorithms for a class of nonlocal models featured by finite domain of nonlocal interactions.
Important applications include the computation of the peridynamics model of solid mechanics and models of nonstandard diffusion.
The main message is that the far field behavior of the nonlocal interaction kernel is decisive in 
order for the fast computation of nonlocal models. 
We use numerical tests for the nonlocal diffusion model as a proof of concept, and demonstrate that better regularity of the kernel away from zero is better 
in the interest of fast computation. This is reasonable since most fast algorithms that can be applied to 
non-translation invariant operators and general meshes 
are based on low rank approximation of the far field interactions. 
Based on this observation, the kernels that are most often used in the peridynamics simulations fail to work under
existing fast algorithms since they all have non-smooth truncations at a finite distance. 
We propose a remedy for such kernels, and that is the splitting of the singularities at origin and at the finite distance and treat them separately. 
The splitting is done by constructing polynomials that match the derivatives of the original kernels at finite distances. The second type of singularity after splitting is on larger sets than just point discontinuity. 
We propose a new FMM like fast algorithm to deal with such singularities which exploits the smoothness of the kernels away from the truncation. 
Such idea that uses low rank approximations of kernels for different regions is also presented in \cite{Engquist2007}
for a completely different problem. 
Note that we only focus on the fast matrix-vector multiplication in this work, and
it can be used in a Krylov subspace method if computing the inverse of a matrix is the interest.
Notice that the kernel function that uses to regularize the original kernel typically is not small,
as indicate in the bottom right picture in Fig.~\ref{fig:1}. If the original kernel is a small perturbation of a kernel that have smooth decay 
in the far field, 
then the inverse of the corresponding matrix could be simply computed using the series expansion under the smallest assumption.
The goal of this paper is to show that fast algorithms for nonlocal operators with finite interaction distances are possible.
Even if the techniques developed here substantially reduces the computational cost, there is room for 
improvements from future research. It would be natural to pursue further research as nonlocal modeling becomes more
common. The essentially optimal computational complexity of $O(N\log N)$ is for 1d
and in higher dimensions there is algebraic complexity higher than linear. 
For further complexity reduction in higher dimensions a better representation 
of precomputed quantities near the horizon would be required and it would 
be natural to use different geometries than boxes, for example, curvelets \cite{Candes2000} in 
2d deals efficiently with line discontinuities.

Based on our findings, several more points have the value of further study. 
First, our work offers important messages to nonlocal modeling from the perspective of scientific computing. 
It is worthwhile to explore in the future whether the kernels with smooth decays at far field can be used in replace of  
non-smooth kernels in the nonlocal models without compromising the physics. 
Second, 
in order to understand the properties 
of the algorithms, we separated the global HIF technique from the handling of 
discontinuous at the horizon by using a simple splitting. The more interesting question on
the optimal way to split the singularities could be asked in the future. 
Moreover, it will be more desirable to build a tighter coupling of the singularities into the original HIF technique, 
and the question on detecting 
the zero values of the kernel could be asked as mentioned earlier.
Finally, we note that our case study is restricted to nonlocal diffusion type problem. 
It is reasonable that for the peridynamics system of equations similar studies can 
be done.

\section*{Acknowledgements}
The first author want to thank Lexing Ying and Yuwei Fan for providing with helpful discussions and HIF codes 
during her visit to Stanford University. 
The authors also thank Qiang Du for helpful discussions on the subject. 
The Oden Institute  is acknowledged for its support.  
The research of Bjorn Engquist is supported in part by the U.S. NSF grant DMS-1620396. The research of Xiaochuan Tian is supported in part
by the U.S. NSF grant DMS-1819233. We thank the referees for the helpful comments to improve this work. 

\bibliographystyle{siam}
\bibliography{txc}

             \end{document}